\documentclass[12pt]{article}
\usepackage{exscale,relsize}
\usepackage{comment}
\usepackage{lscape}
\usepackage{amsmath}
\usepackage[sort,nocompress]{cite}
\usepackage{amsfonts}
\usepackage[hidelinks]{hyperref}
\usepackage{amssymb}
\usepackage{calc}
\usepackage{theorem}
\usepackage{pifont}      
\usepackage{array}
\usepackage{mathtools} 

\usepackage{color}

\usepackage{graphicx} 
\usepackage{float} 
\usepackage{subcaption} 

\newcommand{\sgn}{\ensuremath{\operatorname{sgn}}}
\newcommand{\N}{\mathbb{N}}   

\newcommand{\norm}[1]{\left\lVert#1\right\rVert} 
\newcommand{\ip}[2]{\langle#1,#2\rangle} 

\newcommand{\gph}{\ensuremath{\operatorname{gph}}}

\newcommand{\R}{\ensuremath{\mathbb R}}
\newcommand{\Rn}{\ensuremath{\mathbb R^n}}

\newcommand{\dom}{\ensuremath{\operatorname{dom}}}

\newcommand{\inte}{\ensuremath{\operatorname{int}}}

\newcommand{\ran}{\ensuremath{\operatorname{ran}}}

\newcommand{\Id}{\ensuremath{\operatorname{Id}}}

\newcommand{\prox}{\ensuremath{P}}
\newcommand{\proj}{\ensuremath{\operatorname{Proj}}}

\newtheorem{theorem}{Theorem}[section]

\newtheorem{fact}[theorem]{Fact}
\newtheorem{corollary}[theorem]{Corollary}
\newtheorem{proposition}[theorem]{Proposition}
\newtheorem{defn}[theorem]{Definition}

\theoremstyle{plain}{\theorembodyfont{\rmfamily}
	}
\theoremstyle{plain}{\theorembodyfont{\rmfamily}
	}
\theoremstyle{plain}{\theorembodyfont{\rmfamily}
	}
\theoremstyle{plain}{\theorembodyfont{\rmfamily}
	\newtheorem{example}[theorem]{Example}}
\theoremstyle{plain}{\theorembodyfont{\rmfamily}
	\newtheorem{remark}[theorem]{Remark}}
\theoremstyle{plain}{\theorembodyfont{\rmfamily}
	}
\def\proof{\noindent{\it Proof}. \ignorespaces}
\def\endproof{\ensuremath{\quad \hfill \blacksquare}}

\newcommand{\pluss}{{\hskip1pt \raise1pt\vbox{\hrule width6pt \vskip1pt
			\hrule width6pt}\kern-4pt{\lower1pt\hbox{\vrule height6pt \kern1pt\vrule
				height6pt}}\hskip5pt}}

\newcommand{\argmin}{\mathop{\rm argmin}\limits}

\newcommand{\levp}{\partial_p^\lambda}
\newcommand{\sph}{\mathbb{S}}
\usepackage{pgfplots}
\usepackage{tikz}
\usepackage{subcaption}

\begin{document}
	\title{Every proximal mapping is a resolvent of level proximal subdifferential}
\author{Xianfu Wang\footnote{Department of Mathematics, University of British Columbia, Kelowna, B.C. V1V 1V7, Canada. E-mail: shawn.wang@ubc.ca.} \ and Ziyuan Wang\footnote{Department of Mathematics, University of British Columbia, Kelowna, B.C. V1V 1V7, Canada. E-mail: ziyuan.wang@alumni.ubc.ca.}}

\date{March 2, 2023}
	\maketitle

	\begin{abstract} We propose a level proximal subdifferential for a proper lower semicontinuous function. 
Level proximal subdifferential is a uniform refinement of the well-known proximal subdifferential, and
has the pleasant feature that its resolvent always coincides with the proximal mapping of a function.
It turns out that the resolvent representation of proximal mapping in terms of Mordukhovich limiting subdifferential is only valid for hypoconvex functions.
We also provide properties of level proximal subdifferential and numerous examples to illustrate our results.
	\end{abstract}
{
\noindent
{\bfseries 2020 Mathematics Subject Classification:}
{Primary 49J52, 49J53; Secondary 47H05, 47H04, 26B25.
}

\noindent {\bfseries Keywords:} Proximal mapping, hypoconvex function, level proximal subdifferential,
Mordukhovich limiting subdifferential.
}

	\section{Introduction}
Throughout $\Rn$ denotes the standard $n$-dimensional Euclidean space with inner product $\langle\cdot,\cdot\rangle$ and
norm $\|x\|=\sqrt{\langle x, x\rangle}$ for $x\in\Rn$.
Let $f:\Rn\to\overline\R$ be proper and lower semicontinuous~(lsc). The proximal mapping (or operator) of $f$ with parameter $\lambda>0$ is $$\prox_\lambda f:\Rn\to 2^{\Rn}: x\mapsto\argmin_{y\in\Rn}\bigg\{f(y)+\frac{1}{2\lambda}\norm{y-x}^2\bigg\},$$
which is a fundamental object in continuous optimization and variational analysis. {For important properties and applications of proximal mappings, see \cite{rockafellar_variational_1998,bt05,cwp20} on general functions, \cite{beck2017first,BC} on convex functions, and the references therein.}

When the given function is ``nice", its proximal operator admits a pleasant representation.
\begin{fact}\emph{\cite[Proposition 12.19]{rockafellar_variational_1998}}\label{fact: weak cvx} Let $\lambda>0$ and let $f:\Rn\to\overline{\R}$ be proper, lsc and $1/\lambda$-hypoconvex, that is, $f+(1/2\lambda)\norm{\cdot}^2$ is convex. Then
	\begin{equation}\label{eq: prox identity}
\prox_\lambda f=\left(\Id+\lambda\partial f\right)^{-1},
	\end{equation}
where $\partial f$ denotes the Mordukhovich limiting subdifferential operator of $f$ (see Definition~\ref{Defn: subdifferentials}).
\end{fact}
However identity~(\ref{eq: prox identity}) may fail beyond the usual hypoconvexity setting.
For example, consider the function $f:\R\to\R$ given by $f(0)=0$ and $f(x)=1$ if $x\neq0$, which is the one-dimensional zero ``norm". Then $(\forall |x|<\sqrt{2\lambda})$ $\prox_\lambda f(x)=\{0\}\subset\{0,x\}=\left(\Id+\lambda\partial f\right)^{-1}(x)$; see Example~\ref{eg: zero no} for details. Such an observation reveals that resolvents of limiting subdifferential are in general too big to capture proximal operators, leading to the question below:
\begin{center}
	\emph{Which subdifferential completely represents proximal operator?}
\end{center}
We resolve this question completely in this paper. We begin by showing that the converse of Fact~\ref{fact: weak cvx} holds, suggesting that the limiting subdifferential representation in~(\ref{eq: prox identity}) is only valid within the class of hypoconvex functions; see Theorem~\ref{thm: weak cvx equivalency}.
Then we depart from the literature by introducing the level proximal subdifferential (see Definition~\ref{defn: levp}), a new concept that refines the classical proximal subdifferential. In striking contrast to the literature, the resolvent of level proximal subdifferential always coincides with the associated proximal operator, regardless of the presence of hypoconvexity.
Several useful properties and examples of the level proximal subdifferential are also given to illustrate its advantages.


The paper is organized as follows.
After necessary notation and preliminaries in Section~\ref{sec: notation}, main results are presented in Section~\ref{sec: main}. Then we illustrate our results in Section~\ref{sec: examples} by examples.
\section{Notation and preliminaries}\label{sec: notation}


Our notation is standard and follows \cite{rockafellar_variational_1998}. The extended real line is $\overline{\R}=\R\cup\{\pm\infty\}$.
We say that a function $f:\Rn\to\overline{\R}$ is proper, if $f>-\infty$ and domain $\dom f=\{x\in\Rn: f(x)<\infty\}$ is nonempty.
A proper, lsc function $f:\Rn\to\overline{\R}$ is prox-bounded if there exists $\lambda>0$ such that $f+(1/2\lambda)\norm{\cdot}^2$ is bounded below; see, e.g.,~\cite[Exercise 1.24]{rockafellar_variational_1998}.
The supremum of the set of all such $\lambda$ is the threshold $\lambda_f$ of prox-boundedness for $f$.
The indicator function of a set $K\subseteq\Rn$ is $\delta_K(x)=0$ if $x\in K$ and $\delta_K(x)=\infty$ otherwise.
We write $\prox_{1}\delta_{K}=\proj_{K}$.
Let $\lambda>0$. Then an operator $A:\Rn\to2^{\Rn}$ is $1/\lambda$-hypomonotone if $(\forall (x,u),(y,v)\in\gph A)~\ip{x-y}{u-v}\geq-(1/\lambda)\norm{x-y}^2$; a function $f:\Rn\to\overline{\R}$ is $1/\lambda$-hypoconvex if $f+(1/2\lambda)\norm{\cdot}^2$ is convex.

 We will use frequently the following concepts from variational analysis; see, e.g., \cite{rockafellar_variational_1998, mor2018variational,clarke98}.
 \begin{defn}\label{Defn: subdifferentials}
	Let $f:\Rn\rightarrow\overline{\R}$ be a proper function. We say that
	\begin{enumerate}
	\item $v\in\Rn$ is a proximal subgradient of $f$ at $\bar{x}\in\dom f$, denoted by $v\in\partial_p f(\bar{x})$, if there exist $\rho>0$ and $\delta>0$ such that
	\begin{equation}
	(\forall \norm{x-\bar{x}}\leq\delta)~f(x)\geq f(\bar{x})+\ip{v}{x-\bar{x}}-\frac{\rho}{2}\norm{x-\bar{x}}^2.
	\end{equation}
	\item $v\in\Rn$ is a \textit{Fr\'echet subgradient} of $f$ at $\bar{x}\in\dom f$, denoted by $v\in\hat{\partial}f(\bar{x})$, if
		\begin{equation}\label{Formula:frechet subgradient inequality}
			f(x)\geq f(\bar{x})+\ip{v}{x-\bar{x}}+o(\norm{x-\bar{x}}).
		\end{equation}
	\item  $v\in\Rn$ is a \textit{Mordukhovich (limiting) subgradient} of $f$ at $\bar{x}\in\dom f$, denoted by $v\in\partial f(\bar{x})$, if
	\begin{equation}\label{Formula:limiting subgraident definition}
		v\in\{v\in\Rn:\exists x_k\xrightarrow[]{f}\bar{x},\exists v_k\in\hat{\partial}f(x_k),v_k\rightarrow v\},
	\end{equation}
	where $x_k\xrightarrow[]{f}\bar{x}\Leftrightarrow x_k\rightarrow\bar{x}\text{ and }f(x_k)\rightarrow f(\bar{x})$. Moreover, we set $\dom\partial f=\{x\in\Rn:\partial f(x)\neq\emptyset\}$.
	\end{enumerate}
\end{defn}
Clearly, $\partial_{p}f(x)\subseteq \hat{\partial}f(x)\subseteq\partial f(x)$ for every $x\in\Rn$.

\section{Level proximal subdifferential and its resolvent}\label{sec: main}

In this section, we present main results of this paper, beginning with an equivalence characterization of the identity $\prox_\lambda f  =\left(\Id+\lambda\partial f\right)^{-1}$.  One direction of the following result is known from~\cite[Proposition 12.19]{rockafellar_variational_1998}. However, to the best of our knowledge, the converse direction is new.

\begin{theorem}\label{thm: weak cvx equivalency} Let $f:\Rn\to\overline{\R}$ be proper, lsc and prox-bounded with threshold $\lambda_f\in(0,\infty]$. Then for every $\lambda\in(0,\lambda_f)$
	\begin{equation}\label{equivalence1}
		\prox_\lambda f  =\left(\Id+\lambda\partial f\right)^{-1}\Leftrightarrow f+\frac{1}{2\lambda}\norm{\cdot}^2\text{ is convex. }
	\end{equation}
\end{theorem}

\proof Suppose that $\prox_\lambda f=\left(\Id+\lambda\partial f\right)^{-1}$. Then appealing to~\cite[Proposition 12.19]{rockafellar_variational_1998} yields that $\left(\Id+\lambda\partial f\right)^{-1}$ is monotone, and so is the set-valued operator $\Id+\lambda\partial f$ by~\cite[Exercise 12.4]{rockafellar_variational_1998}. In turn $\partial\left((1/2\lambda)\norm{\cdot}^2+f \right)=(1/\lambda)\Id+\partial f$ is monotone. Invoking~\cite[Theorem 12.17]{rockafellar_variational_1998} implies the right hand side of~(\ref{equivalence1}). The converse holds due to~\cite[Proposition 12.19]{rockafellar_variational_1998}.~\endproof

Having established equivalence~(\ref{equivalence1}),  we now study proximal operator representation beyond the setting of Theorem~\ref{thm: weak cvx equivalency}. To this end, we propose a new subdifferential, which is a refinement of $\partial_{p}f$.


\begin{defn}[level proximal subdifferential] \label{defn: levp} Let $f:\Rn\to\overline{\R}$ be proper, $x\in\dom f$ and let $\lambda>0$. We say that $u\in\Rn$ is a $\lambda$-level proximal subgradient of $f$ at $x$, denoted by $u\in\levp f(x)$, if
	\begin{equation}\label{lep ineq}
		(\forall y\in\Rn)~f(y)\geq f(x)+\ip{u}{y-x}-\frac{1}{2\lambda}\norm{y-x}^2.
 	\end{equation}

\end{defn}
\begin{remark}\label{rem: compare with known}
	\begin{enumerate}
		\item
		{
			Inequality~(\ref{lep ineq}) with a sufficiently small $\lambda>0$ has appeared in the proof of \cite[Proposition 8.46(f)]{rockafellar_variational_1998} as a technical necessity to globalize proximal subgradient inequality.
		   However, level proximal subdifferential, which emphasizes on inequality~(\ref{lep ineq}) with a predetermined parameter, was not investigated.
	}
		\item
		{
			We mention another connection between level proximal subdifferential and the literature.
			Let $\lambda>0$ and let $U\subset\Rn$ be a bounded open convex set. In~\cite[Theorem 5.1]{clarke_proximal_1995}, Clarke, Stern and Wolenski investigated the following inequality: for each $x\in U$
		\[(\exists u\in\Rn)
		(\forall y\in U)~
		f(y)\geq f(x)
		+
		\ip{u}{y-x}
		-
		\frac{1}{2\lambda}\norm{y-x}^2,
		\]
which clearly holds for $u\in\levp f(x)$.
	}
	\end{enumerate}
\end{remark}
Simple algebraic manipulations of \eqref{lep ineq} give the following characterizations of $\partial_{p}^{\lambda}f$.
\begin{proposition}\label{p:basic} Let $u,x\in\Rn$. The following are equivalent:
\begin{enumerate}
\item $u\in \partial_{p}^{\lambda}f(x).$
\item $x\in \prox_{\lambda}f(x+\lambda u),$ i.e.,
$$(\forall y\in\Rn)\ f(y)+\frac{1}{2\lambda}\|(x+\lambda u)-y\|^2\geq f(x)+\frac{1}{2\lambda}\|(x+\lambda u)-x\|^2.$$
\item\label{i:fenchel}
 $u+x/\lambda$ is a Fenchel subgradient \cite{bh96,beck2017first} of $\displaystyle f+\frac{1}{2\lambda}\|\cdot\|^2$ at $x$, i.e.,
$$(\forall y\in\Rn)\  f(y)+\frac{1}{2\lambda}\|y\|^2\geq f(x)+\frac{1}{2\lambda}\|x\|^2+
\left\langle u+\frac{x}{\lambda}, y-x \right\rangle.$$
\end{enumerate}
\end{proposition}
It is well-know that if a proper lsc function $f:\Rn\rightarrow \overline{\R}$ is
Fenchel subdifferentiable on a convex subset $C$ of $\Rn$, then
$f$ is convex on $C$; see, e.g., \cite[Lemma 3.11]{beck2017first}.

\begin{proposition}\label{p:hypoconvex} On every nonempty convex subset of $\dom\partial_{p}^{\lambda}f$, the function
$f$ is $1/\lambda$-hypoconvex. In particular, if
$\dom\partial_{p}^{\lambda}f=\dom f$ and $\dom f$ is convex, then $f$ is $1/\lambda$-hypoconvex.
\end{proposition}
\begin{proof} Let $C$ be a nonempty convex subset of $\dom\partial_{p}^{\lambda}f$. At every $x\in C$,
by Proposition~\ref{p:basic}\ref{i:fenchel}, the function $\displaystyle f+(1/2\lambda)\|\cdot\|^2$
is Fenchel
subdifferentiable at $x$. Since $C$ is convex,
we conclude that $\displaystyle f+(1/2\lambda)\|\cdot\|^2$ is convex on $C$. The remaining result
is immediate.
\end{proof}

Several useful properties of $\levp f$ are in order.
\begin{proposition}\label{prop: basic propeties}Let $f:\Rn\to\overline{\R}$ be proper and lsc. Then the following hold:
\begin{enumerate}
	\item\label{convex and closed} $(\forall x\in\dom f)~(\forall \lambda>0)$ $\partial_p^{\lambda}f(x)$ is convex and closed.
	\item\label{relation1} $(\forall 0<\lambda_1\leq\lambda_2)~(\forall x\in\dom f)$ $\partial_p^{\lambda_2}f(x)\subseteq\partial_p^{\lambda_1}f(x)$ and $\partial_p f(x)=\bigcup_{\lambda>0}\levp f(x)$.
	\item\label{range and domain} If $(\exists\lambda>0)$ $\ran\levp f$ is bounded, then $\dom\levp f$ is closed.
	\item\label{relation2} $(\forall x\in\dom f)~(\forall \lambda>0)$ $\levp f(x)\subseteq \partial_p f(x)\subseteq\partial f(x)$.
	\item\label{convex subdifferential} $(\forall x\in\dom f)$ $\bigcap_{\lambda>0}\levp f(x)=\{v\in\Rn:(\forall y\in\Rn)~f(y)\geq f(x)+\ip{v}{y-x}\}$, the Fenchel subdifferential of $f$ at $x$.
	\item\label{gradient} If $f$ is differentiable on $\inte\dom f\neq\emptyset$ and $1/\lambda$-hypoconvex, then $(\forall x\in\inte\dom f)$ $\levp f(x)=\{\nabla f(x)\}$.
\end{enumerate}
\end{proposition}

\proof \ref{convex and closed} Immediate from the definition.

\ref{relation1} Pick $(x,u)\in\gph\partial_p^{\lambda_2}f$. Then $(\forall y\in\Rn)~f(y)\geq f(x)+\ip{u}{y-x}-(1/2\lambda_2)\norm{y-x}^2\geq f(x)+\ip{u}{y-x}-(1/2\lambda_1)\norm{y-x}^2$, implying that $\partial_p^{\lambda_2}f(x)\subseteq\partial_p^{\lambda_1}f(x)$. Invoking~\cite[Proposition 8.46]{rockafellar_variational_1998} justifies that $\partial_p f(x)=\bigcup_{\lambda>0}\levp f(x)$.

\ref{range and domain} Let $(x_k,u_k)_{k\in\N}\subseteq\gph\levp f$ be such that $x_k\to x\in\Rn$. Taking subsequence if necessary, assume without loss of generality that $u_k\to u$. Then
\begin{align*}
(\forall y\in\Rn)~f(y)&\geq\liminf_{k\to\infty}\left( f(x_k)+\ip{u_k}{y-x_k}-\frac{1}{2\lambda}\norm{y-x_k}^2\right)\\
&\geq f(x)+\ip{u}{y-x}-\frac{1}{2\lambda}\norm{y-x}^2,
\end{align*}
where the last inequality holds because $f$ is lsc, entailing the desired claim.

\ref{relation2} The first inclusion follows from \ref{relation1}. Then apply~\cite[Proposition 8.46(e)]{rockafellar_variational_1998}.

\ref{convex subdifferential} For simplicity, denote $A=\{v\in\Rn:(\forall y\in\Rn)~f(y)\geq f(x)+\ip{v}{y-x}\}$. First we show that $\bigcap_{\lambda>0}\levp f(x)\subseteq A $. Assume without loss of generality that  $\bigcap_{\lambda>0}\levp f(x)$ is nonempty, otherwise it is trivial. Let $u\in\bigcap_{\lambda>0}\levp f(x)$. Then $$(\forall \lambda>0)~(\forall y\in\Rn)~f(y)\geq f(x)+\ip{u}{y-x}-\frac{1}{2\lambda}\norm{y-x}^2\Rightarrow f(y)\geq f(x)+\ip{u}{y-x},$$
where the implication holds by taking $\lambda\to\infty$. Next we justify $ A \subseteq\bigcap_{\lambda>0}\levp f(x)$. Assume without loss of generality again that $ A $ is nonempty and pick $u\in A $. Then $(\forall\lambda>0)$ $(\forall y\in\Rn)~f(y)\geq f(x)+\ip{u}{y-x}\geq f(x)+\ip{u}{y-x}-(1/2\lambda)\norm{y-x}^2$, completing the proof.

\ref{gradient} Note that $\levp f(x)\subseteq\partial_p f(x)\subseteq\{\nabla f(x)\}$. Thus it suffices to show that $\nabla f(x)\in\levp f(x)$. By assumption,
\[(\forall y\in\Rn)~f(y)+\frac{1}{2\lambda}\norm{y}^2\geq f(x)+\frac{1}{2\lambda}\norm{x}^2+\ip{\nabla f(x)+\frac{1}{\lambda}x}{y-x},\]
implying~(\ref{lep ineq}) after rearrangement.
\endproof

Amazingly, the resolvent of level proximal subdifferential always coincides with proximal mapping, regardless of hypoconvexity. Moreover, the graph of level proximal subdifferential can be obtained by linearly transforming the graph of proximal operator, and vice versa.
\begin{theorem}[resolvent representation of $\prox_\lambda f$]\label{thm: prox identity}
Let $f:\Rn\to\overline{\R}$ be proper, lsc and prox-bounded with threshold $\lambda_f\in(0,\infty]$. Then
\begin{equation}\label{prox and prox sub}
(\forall 0<\lambda<\lambda_f)~\prox_\lambda f =\left(\Id+\lambda\levp f\right)^{-1}.
\end{equation}
Consequently, operator $\left(\Id+\lambda\levp f\right)^{-1}:\Rn\to2^{\Rn}$ has full domain, is compact-valued, and for every $x\in\Rn$
\begin{align*}
	\left(\Id+\lambda\partial_p^\lambda f\right)^{-1}(x)=\big\{v\in\Rn:v_k\to v,\text{where }(\exists x_k\to x)(\exists \lambda_k\to\lambda)v_k\in\left(\Id+\lambda_k\partial_p^{\lambda_k} f\right)^{-1}(x_k)\big\}.
\end{align*}
Furthermore $\gph\levp f=T\left(\gph \prox_\lambda f\right)$, where
$T=\begin{bmatrix}
	0&\Id\\\Id/\lambda&-\Id/\lambda
\end{bmatrix}$ is an invertible linear transformation.
\end{theorem}
\proof Note that $(\forall x\in\Rn)$~$\prox_\lambda(x)$ is nonempty thanks to~\cite[Theorem 1.25]{rockafellar_variational_1998}.  By the definition of proximal operator,
\begin{align*}
(\forall x\in\Rn)~x^*\in\prox_\lambda f (x)&\Leftrightarrow( \forall y\in\Rn)~f(y)+\frac{1}{2\lambda}\norm{y-x}^2\geq f(x^*)+\frac{1}{2\lambda}\norm{x^*-x}^2\\
&\Leftrightarrow( \forall y\in\Rn)~f(y)\geq f(x^*)+\frac{1}{\lambda}\ip{x-x^*}{y-x^*}-\frac{1}{2\lambda}\norm{y-x^*}^2\\
&\Leftrightarrow \frac{1}{\lambda}(x-x^*)\in\levp f(x^*)\Leftrightarrow x^*\in\left(\Id+\lambda\levp f\right)^{-1}(x),
\end{align*}
from which~(\ref{prox and prox sub}) and the liner transformation result hold. Apply~\cite[Theorem 1.25]{rockafellar_variational_1998} again to complete the proof.~\endproof

\begin{corollary} Let $f:\Rn\to\overline{\R}$ be proper, lsc and prox-bounded with threshold $\lambda_f\in(0,\infty]$.
Then $\partial_p^{\lambda}f$ is $1/\lambda$-hypomonotone on its domain $\dom \partial_p^{\lambda}f$.
\end{corollary}
\begin{proof}
Apply equation~\eqref{prox and prox sub} of Theorem~\ref{thm: prox identity}.
\end{proof}

A systematic study on properties of resolvents of hypomontone operators has been given in \cite{bmw21}.
{In the presence of hypoconvexity and convexity, level proximal subdifferential is indistinguishable from other known subdifferentials.}

\begin{theorem}[subdifferential relationships] \label{thm: relations}
	Let $f:\Rn\to\overline{\R}$ be proper, lsc and prox-bounded with threshold $\lambda_f\in(0,\infty]$. Then the following hold:
	\begin{enumerate}
		\item\label{hypoconvex} Let $0<\lambda<\lambda_f$. Then $f$ is $1/\lambda$-hypoconvex $\Leftrightarrow\levp f=\partial f$ $\Leftrightarrow\levp f=\partial_p f\Leftrightarrow\partial_p^\lambda f=\hat\partial f$. 
		\item\label{convex} Suppose that $\lambda_f=\infty$. Then $f$ is convex $\Leftrightarrow (\forall \lambda>0)$ $\levp f=\partial f$ $\Leftrightarrow (\forall \lambda>0)$ $\levp f=\partial_p f$  $\Leftrightarrow (\forall \lambda>0)$ $\levp f=\hat\partial f$.
	\end{enumerate}
\end{theorem}
\proof \ref{hypoconvex} 
Assume that $\levp f=\partial f$. Then Theorem~\ref{thm: prox identity} yields that $\prox_\lambda f =\left(\Id+\lambda\partial f\right)^{-1}$, which is monotone by using~\cite[Proposition 12.19]{rockafellar_variational_1998}. In turn $\Id+\lambda \partial f$ is monotone by~\cite[Exercise 12.4]{rockafellar_variational_1998} and consequently $\partial f$ is $1/\lambda$-hypomonotone. The hypoconvexity of $f$ follows immediately; see, e.g.,~\cite[Exercise 12.61]{rockafellar_variational_1998}.
To justify the converse direction, observe that $$\left(\Id+\lambda\levp f\right)^{-1}=\prox_\lambda f=\left(\Id+\lambda\partial f\right)^{-1}$$ by Theorems~\ref{thm: weak cvx equivalency} and~\ref{thm: prox identity}, implying that $\partial_p^\lambda f=\partial f$.

We now show that $\levp f=\partial_p f$ implies $f$ being $1/\lambda$-hypoconvex. Note that $\left(\Id+\lambda\partial_p f\right)^{-1}$ is monotone due to Theorem~\ref{thm: prox identity} and~\cite[Exercise 12.4]{rockafellar_variational_1998}, entailing $\partial_p f$ to be $1/\lambda$-hypomonotone. Take arbitrary $(x,u),(y,v)\in\gph\partial f$. Then there exist $(x_k,u_k),(y_k,v_k)\in\gph\partial_p f$ such that $(x_k,u_k)\to (x,u)$ and $(y_k,v_k)\to(y,v)$. In turn,
\begin{equation*}
(\forall k\in\N)~\ip{x_k-y_k}{u_k-v_k}\geq-\frac{1}{\lambda}\norm{x_k-y_k}^2\Rightarrow \ip{x-y}{u-v}\geq-\frac{1}{\lambda}\norm{x-y}^2,
\end{equation*}
which justifies the desired claim by invoking~\cite[Exercise 12.61]{rockafellar_variational_1998} again. The converse direction is immediate by using Proposition~\ref{prop: basic propeties}\ref{relation2}. The last equivalence regarding $\hat\partial f$ can be proved by a similar argument.

\ref{convex} Suppose that $f$ is convex. Then $(\forall \lambda >0)$ $\partial f(x)\subseteq\partial_p^\lambda f(x)\subseteq\partial _pf(x)\subseteq\partial f(x)$, where the first inclusion is implied by~(\ref{lep ineq}) and the second one holds due to Proposition~\ref{prop: basic propeties}, from which the desired result readily follows. Now suppose that $(\forall \lambda>0)$ $\partial_p^\lambda f=\partial f$. Then~\ref{hypoconvex} implies that for every $(x,u),(y,v)\in\gph\partial f$
\begin{align*}
(\forall \lambda>0)~\ip{x-y}{u-v}\geq-\frac{1}{\lambda}\norm{x-y}^2\Rightarrow\ip{x-y}{u-v}\geq0,
\end{align*}
where the right hand side holds by taking $\lambda\to\infty$, which means that $\partial f$ is monotone. It follows immediately that $f$ is convex; see, e.g.,~\cite[Theorem 12.17]{rockafellar_variational_1998}. The remainder can be justified by a similar argument as in statement\ref{hypoconvex}. \endproof

We end this section by investigating sum rules of level proximal subdifferential.

\begin{proposition}[addition of functions] Let $f,g:\Rn\to\overline{\R}$ be proper, $\lambda_1,\lambda_2>0$ and let $\lambda_3=(\lambda_1\lambda_2)/(\lambda_1+\lambda_2)$. Let $x\in\dom(f+g)$. Then the following hold:
\begin{enumerate}
	\item\label{weak sum rule} $\partial_p^{\lambda_1}f(x)+\partial_p^{\lambda_2}g(x)\subseteq\partial_p^{\lambda_3}(f+g)(x)$.
	\item\label{strong sum rule}  Suppose that $\nabla g$ is $1/\lambda_2$-Lipschitz. Then
\begin{equation}\label{sum rule other side}
	\partial_p^{\lambda_3}(f+g)(x)\subseteq \partial_p^{\hat\lambda_3}f(x)+\nabla g(x)=\partial_p^{\hat\lambda_3}f(x)+\partial_p^{\lambda_2} g(x),
\end{equation}
	where $\hat\lambda_3=(\lambda_2\lambda_3)/(\lambda_2+\lambda_3)$.
Assume in addition that $f$ is $1/\lambda_1$-hypoconvex. Then
	 \begin{equation}\label{sum rule equality}
		\partial_p^{\lambda_1}f(x)+\partial_p^{\lambda_2} g(x)=\partial_p^{\lambda_1}f(x)+\nabla g(x)=\partial_p^{\lambda_3}(f+g)(x).
	\end{equation}
\end{enumerate}

\end{proposition}
\proof \ref{weak sum rule} Assume without loss of generality that $\partial_p^{\lambda_1}f(x),\partial_p^{\lambda_2}g(x)$ are nonempty, otherwise the desired inclusion is trivial. Now, pick $u\in\partial_p^{\lambda_1}f(x),v\in\partial_p^{\lambda_2}g(x)$. Then adding corresponding level proximal subgradient inequality~(\ref{lep ineq}) yields that
\begin{align*}
(\forall y\in\Rn)~(f+g)(y)\geq (f+g)(x)+\ip{u+v}{y-x}-\frac{1}{2\lambda_3}\norm{y-x}^2,
\end{align*}
meaning that $u+v\in\partial_p^{\lambda_3}(x)$.

\ref{strong sum rule} Note that the Lipschitz assumption of $\nabla g$ amounts to $(1/2\lambda_2)\norm{\cdot}^2\pm g$ being convex; see, e.g., {~\cite[Lemma 2.5]{wang2022mirror}} or \cite[Lemma 5.7, p. 109]{beck2017first}.
Then Proposition~\ref{prop: basic propeties}\ref{gradient} implies that $\partial_p^{\lambda_2}g(x)=\nabla g(x)$ and $\partial_p^{\lambda_2}(-g)(x)=-\nabla g(x)$.
Pick an arbitrary $u\in\partial_p^{\lambda_3}(f+g)(x)$. Then for every $y\in\Rn$
\begin{align*}
	f(y)&\geq f(x)+\ip{u-\nabla g(x)}{y-x}-\frac{1}{2\lambda_3}\norm{y-x}^2-\left(g(y)-g(x)-\ip{\nabla g(x)}{y-x} \right)\\
	&\geq  f(x)+\ip{u-\nabla g(x)}{y-x}-\left(\frac{1}{2\lambda_3}+\frac{1}{2\lambda_2}\right)\norm{y-x}^2\\
	&=f(x)+\ip{u-\nabla g(x)}{y-x}-\frac{1}{2\hat\lambda_3}\norm{y-x}^2,
\end{align*}
where the second inequality holds because $\partial_p^{\lambda_2}(-g)(x)=\{-\nabla g(x)\}$, which combined with $\partial_p^{\lambda_2}g(x)=\{\nabla g(x)\}$ completes the proof of~(\ref{sum rule other side}).
Finally, we turn to~(\ref{sum rule equality}). Note that $\lambda_1>\hat\lambda_3$. Invoking Theorem~\ref{thm: weak cvx equivalency} and Proposition~\ref{prop: basic propeties}\ref{relation1}\&\ref{relation2}, one concludes that
\[
\partial_p f(x)=\partial_p^{\lambda_1}f(x)\subseteq\partial_p^{\hat\lambda_3}f(x)\subseteq\partial_p f(x)\Rightarrow\partial_p^{\lambda_1}f(x)=\partial_p f(x)=\partial_p^{\hat\lambda_3}f(x),
\]
which entails~(\ref{sum rule equality}) by taking statement\ref{weak sum rule} and~(\ref{sum rule other side}) into account.
~\endproof

\section{Examples}\label{sec: examples}

In this final section, we illustrate our main results by calculating resolvents of level proximal, proximal and limiting subdifferentials of one-dimensional and multi-variable functions. They highlight the differences among these resolvents.

\begin{example}\label{eg: indicator} Let $\lambda>0$ and $f(x)=\delta_{\{-1,1\}}(x)$. Then the following hold:
	\begin{enumerate}
		\item $\partial_p f(x)=\partial f(x)=\begin{cases}
			\R, &\text{if }x\in\{-1,1\};\\
			\emptyset, &\text{otherwise},
		\end{cases}$ and $\levp f(x)=
	\begin{cases}
		[-1/\lambda,\infty),&\text{if }x=1;\\
		(-\infty,1/\lambda], &\text{if }x=-1;\\
		\emptyset,&\text{otherwise. }
	\end{cases}$
\item $\prox_\lambda f=\left(\Id+\lambda\partial_p^\lambda f\right)^{-1}\subset \left(\Id+\lambda\partial_p f\right)^{-1}=\left(\Id+\lambda\partial f\right)^{-1}$, where
$(\forall x\in\R) \left(\Id+\lambda\partial f\right)^{-1}(x)=\{-1,1\}$ and
\[\prox_\lambda f(x)=
\begin{cases}
\{-1,1\},&\text{if }x=0;\\
\sgn(x),&\text{if }x\neq0.
\end{cases}
\]

	\end{enumerate}
\begin{figure}[h!]
	\begin{subfigure}[b]{0.45\textwidth}
		\centering
		\resizebox{\linewidth}{!}{
			\begin{tikzpicture}
				\begin{axis}[xmin = -2, xmax = 2,ymin = -2, ymax = 2,grid = both,
					legend style={
						legend columns=2,
						anchor=south,
						at={(0.5,1)}
					},
					xtick={-2,-1,0,1,2},
					ytick={-2,-1,0,1,2},
					xticklabels={$ $,-1,0,1,2}]
					\addplot[domain = -2:2,samples = 10,smooth,thick,red] {1};
					\addplot[domain = -2:2,samples = 10,smooth,thick,red] {-1};
				\end{axis}
		\end{tikzpicture}}
		\caption{$\left(\Id+\lambda\partial f\right)^{-1}$}
	\end{subfigure}
	\begin{subfigure}[b]{0.45\textwidth}
		\centering
		\resizebox{\linewidth}{!}{
			\begin{tikzpicture}
				\begin{axis}[xmin = -2, xmax = 2,ymin = -2, ymax = 2,grid = both,
					legend style={
						legend columns=2,
						anchor=south,
						at={(0.5,1)}
					},
					xtick={-2,-1,0,1,2},
					ytick={-2,-1,0,1,2},
					xticklabels={$ $,-1,0,1,2}]
					\addplot[domain = -2:0,samples = 10,smooth,thick,blue] {-1};
					\addplot[domain =0:2,samples = 10,smooth,thick,blue] {1};
				\end{axis}
			\end{tikzpicture}
		}
		\caption{$\left(\Id+\lambda\levp f\right)^{-1}$}
	\end{subfigure}
	\caption{Resovlents in Example~\ref{eg: indicator}.}
\end{figure}
\end{example}
\proof (i)  Here we only justify the formula for $\levp f$. Let $v\in\R$. Then $v\in\levp f(1)\Leftrightarrow f (-1)\geq f(1)+2v-2/\lambda\Leftrightarrow v\geq -1/\lambda$. Similarly, $v\in\levp f(-1)\Leftrightarrow f(1)\geq f(-1)+2v-2/\lambda\Leftrightarrow v\leq1/\lambda$.

(ii) Resolvent computation is immediate, hence omitted.
Observe that $\prox_\lambda f$ is the projector onto $\{-1,1\}$ to see the proximal operator formula.
\endproof

\begin{example}\label{eg: zero no} Let $\lambda>0$ and $f(x)=1$ for $x\neq0$ and $f(0)=0$.
Then the following hold:
\begin{enumerate}
	\item $\partial f(x)=\partial_pf(x)=\begin{cases}
		\R,&\text{if $x=0$};\\
		0,&\text{if $x\neq0$,}
	\end{cases}$ and $$
\levp f(x)=\begin{cases}
		\emptyset, &\text{if $0<|x|<\sqrt{2\lambda}$};\\
		[-\sqrt{2/\lambda},\sqrt{2/\lambda}],&\text{if $x=0$};\\
		0, &\text{if $|x|\geq\sqrt{2\lambda}$.}\end{cases}$$
	\item 
$\prox_\lambda f =\left(\Id+\lambda\levp f\right)^{-1}\subset\left(\Id+\lambda\partial_pf\right)^{-1}=\left(\Id+\lambda\partial f\right)^{-1}$, where
$(\forall x\in\R)$ $\left(\Id+\lambda\partial f\right)^{-1}(x)=\left(\Id+\lambda\partial_p f\right)^{-1}(x)=\{0,x\},
$
whereas
	 \[\prox_\lambda f (x)=
	\begin{cases}
		0,&\text{if $|x|<\sqrt{2\lambda}$};\\
		\{0,x\}, &\text{if $|x|=\sqrt{2\lambda}$};\\
	    x,&\text{if $|x|>\sqrt{2\lambda}$}.
	\end{cases}\]
\end{enumerate}
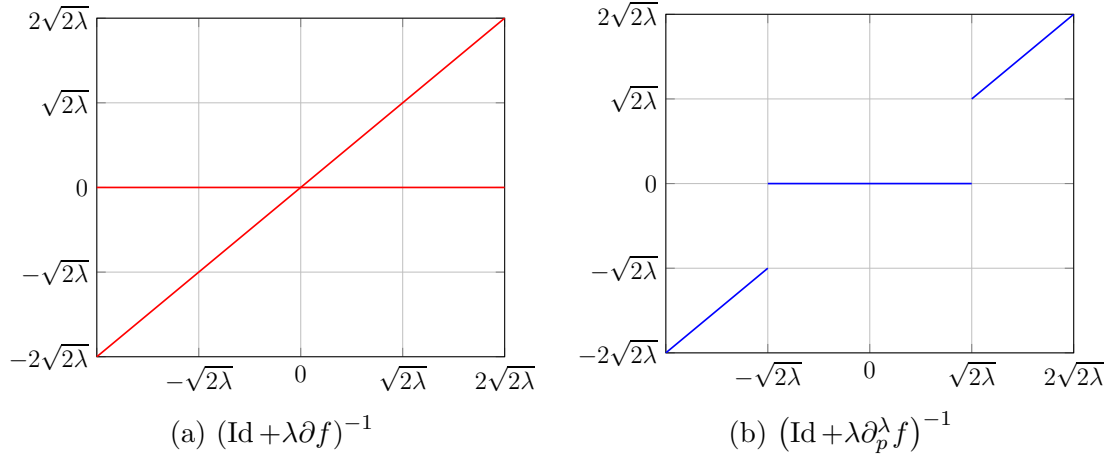
\begin{figure}[h!]
	\begin{subfigure}[b]{0.45\textwidth}
		\centering
		\resizebox{\linewidth}{!}{
			\begin{tikzpicture}
				\begin{axis}[xmin = -2, xmax = 2,ymin = -2, ymax = 2,grid = both,
					legend style={
						legend columns=2,
						anchor=south,
						at={(0.5,1)}
					},
					xtick={-2,-1,0,1,2},
					ytick={-2,-1,0,1,2},
					xticklabels={, $-\sqrt{2\lambda}$,0,$\sqrt{2\lambda}$, $2\sqrt{2\lambda}$},
					yticklabels={$-2\sqrt{2\lambda}$, $-\sqrt{2\lambda}$,0,$\sqrt{2\lambda}$, $2\sqrt{2\lambda}$}]
					\addplot[domain =-2:2,samples = 10,smooth,thick,red] {0};
					\addplot[domain =-2:2,samples = 10,smooth,thick,red] {x};
				\end{axis}
			\end{tikzpicture}
		}
		\caption{$\left(\Id+\lambda\partial f\right)^{-1}$}
	\end{subfigure}
	\begin{subfigure}[b]{0.45\textwidth}
		\centering
		\resizebox{\linewidth}{!}{
			\begin{tikzpicture}
				\begin{axis}[xmin = -2, xmax = 2,ymin = -2, ymax = 2,grid = both,
					legend style={
						legend columns=2,
						anchor=south,
						at={(0.5,1)}
					},
					xtick={-2,-1,0,1,2},
					ytick={-2,-1,0,1,2},
					xticklabels={, $-\sqrt{2\lambda}$,0,$\sqrt{2\lambda}$, $2\sqrt{2\lambda}$},
					yticklabels={$-2\sqrt{2\lambda}$, $-\sqrt{2\lambda}$,0,$\sqrt{2\lambda}$, $2\sqrt{2\lambda}$}]
					
					\addplot[domain =-2:-1,samples = 10,smooth,thick,blue]{x};
					\addplot[domain =1:2,samples = 10,smooth,thick,blue] {x};
					\addplot[domain =-1:1,samples = 10,smooth,thick,blue]{0};
					
				\end{axis}
				
			\end{tikzpicture}
		}
		\caption{$\left(\Id+\lambda\levp f\right)^{-1}$}
	\end{subfigure}
	\caption{Resolvents in Example~\ref{eg: zero no}.}
\end{figure}

\end{example}
\proof (i) Subdifferentials $\partial f$ and $\partial_p f$ follow easily from definition. Note that $v\in\partial_p^\lambda f(0)\Leftrightarrow(\forall y\in\R)~(1/2\lambda)y^2-vy+1\geq0$. Then the discriminant  $\Delta=v^2-4(1/2\lambda)\leq0\Leftrightarrow|v|\leq\sqrt{2/\lambda}$, entailing $\levp f(0)=	[-\sqrt{2/\lambda},\sqrt{2/\lambda}]$.

For $x\neq0$ and $v\in\R$, define $(\forall y\in\R)$ $g(y)=g_{v,x}(y)=f(y)-v(y-x)+(1/2\lambda)(x-y)^2-f(x)$. Then $v\in\levp f(x)\Leftrightarrow (\forall y\in\R)~g(y)\geq0$.  Let $v=0$. Then $(\forall y\in\R)~g(y)=f(y)+(1/2\lambda)(x-y)^2-1$. In turn $(\forall y\in\R)$ $g(y)\geq0\Leftrightarrow (1/2\lambda)x^2-1=g(0)\geq0\Leftrightarrow |x|\geq\sqrt{2\lambda}$, which together with Proposition~\ref{prop: basic propeties}\ref{relation2}  justified the desired formula.

(ii) Note that the proximal formula is one-dimensional hard thresholding operator; see, e.g.,~\cite[Example 6.10]{beck2017first}. The rest follows from simple resolvent computation.
 \endproof

\begin{example}\label{eg: last } Let $\lambda>0$ and $f(x)=0$ if $x\leq0$ and $f(x)=1$ otherwise. Then the following hold:
	\begin{enumerate}
		\item\label{subdifferentials3} $\partial f(x)=\partial_p f(x)=\begin{cases}
			0,&\text{if $x\neq0$};\\
			[0,+\infty),&\text{if $x=0$,}
		\end{cases}$ and $$
\levp f(x)=\begin{cases}
		[0,\sqrt{2/\lambda}],&\text{if $x=0$};\\
		\emptyset,&\text{if $0<x<\sqrt{2\lambda}$};\\
		0,&\text{otherwise.}
	\end{cases}$$
\item 
$\prox_\lambda f =\left(\Id+\lambda\levp f\right)^{-1}\subset\left(\Id+\lambda\partial f\right)^{-1}$, where
\[(\Id+\lambda\partial f)^{-1}(x)=\left(\Id+\lambda\partial_p f\right)^{-1}(x)=
\begin{cases}
\{x,0\},&\text{ if } x\geq0;\\
x,&\text{ if }x<0,
 \end{cases} \]
whereas
 \[\prox_\lambda f (x)=
\begin{cases}
x, &\text{if $x\leq0$ or $x>\sqrt{2\lambda}$};\\
0,&\text{if $0<x<\sqrt{2\lambda}$};\\
\{0,\sqrt{2\lambda}\},&\text{if $x=\sqrt{2\lambda}$}.
\end{cases}\]
	\end{enumerate}

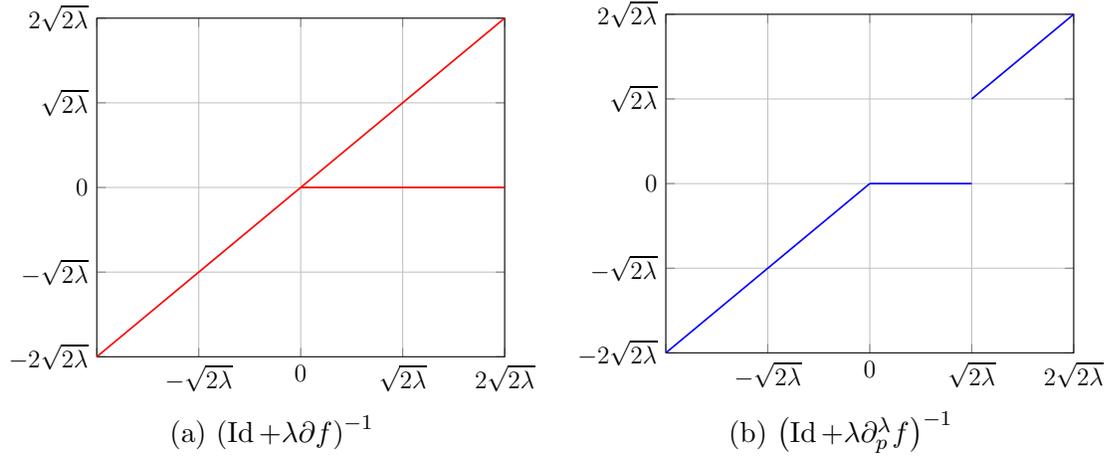
\begin{figure}[h!]
	\begin{subfigure}[b]{0.45\textwidth}
		\centering
		\resizebox{\linewidth}{!}{
			\begin{tikzpicture}
				\begin{axis}[xmin = -2, xmax = 2,ymin = -2, ymax = 2,grid = both,
					legend style={
						legend columns=2,
						anchor=south,
						at={(0.5,1)}
					},
					xtick={-2,-1,0,1,2},
					ytick={-2,-1,0,1,2},
					xticklabels={, $-\sqrt{2\lambda}$,0,$\sqrt{2\lambda}$, $2\sqrt{2\lambda}$},
					yticklabels={$-2\sqrt{2\lambda}$, $-\sqrt{2\lambda}$,0,$\sqrt{2\lambda}$, $2\sqrt{2\lambda}$}]
					\addplot[domain =-2:0,samples = 10,smooth,thick,red] {x};
					\addplot[domain =0:2,samples = 10,smooth,thick,red] {x};
					\addplot[domain =0:2,samples = 10,smooth,thick,red] {0};
				\end{axis}
			\end{tikzpicture}
		}
		\caption{$\left(\Id+\lambda\partial f\right)^{-1}$}
	\end{subfigure}
	\begin{subfigure}[b]{0.45\textwidth}
		\centering
		\resizebox{\linewidth}{!}{
			\begin{tikzpicture}
				\begin{axis}[xmin = -2, xmax = 2,ymin = -2, ymax = 2,grid = both,
					legend style={
						legend columns=2,
						anchor=south,
						at={(0.5,1)}
					},
					xtick={-2,-1,0,1,2},
					ytick={-2,-1,0,1,2},
					xticklabels={, $-\sqrt{2\lambda}$,0,$\sqrt{2\lambda}$, $2\sqrt{2\lambda}$},
					yticklabels={$-2\sqrt{2\lambda}$, $-\sqrt{2\lambda}$,0,$\sqrt{2\lambda}$, $2\sqrt{2\lambda}$}]
					
					\addplot[domain =-2:-0,samples = 10,smooth,thick,blue]{x};
					\addplot[domain =1:2,samples = 10,smooth,thick,blue] {x};
					\addplot[domain =0:1,samples = 10,smooth,thick,blue]{0};
					
				\end{axis}
				
			\end{tikzpicture}
		}
		\caption{$\left(\Id+\lambda\levp f\right)^{-1}$}
	\end{subfigure}
	\caption{Resolvents in Example~\ref{eg: last }.}
\end{figure}

\end{example}
\proof (i) Notice that $v\in\levp f(x)\Leftrightarrow (\forall y\in\R) g(y)\geq0$, where
\[g(y)=
\begin{cases}
	\frac{1}{2\lambda}(y-x)^2-v(y-x)-f(x),&\text{if }y\leq0;\\
	\frac{1}{2\lambda}(y-x)^2-v(y-x)-f(x)+1,&\text{if }y>0.
\end{cases}\]
Consider $x\neq0$, in which case one only needs to check whether $0\in\levp f(x)$ because Proposition~\ref{prop: basic propeties}\ref{relation2}. If $x>0$, then $0\in\levp f(x)$ amounts to $ \min_{y>0}g(y)=1-f(x)\geq0, \min_{y\leq0}g(y)=(1/2\lambda)x^2-1\geq0\Leftrightarrow x\geq\sqrt{2\lambda}$, implying that $\levp f(x)=\emptyset$ when $0<x<\sqrt{2\lambda}$; and $\levp f(x)=\{0\}$ when $x>\sqrt{2\lambda}$. If $x<0$, then $\min_{y\leq0}g(y)=g(x)=-f(x)=0$ and $\inf_{y>0}g(y)=\lim_{y\to0^+}g(y)=(1/2\lambda)x^2+1>0$, so $0\in\levp f(x)$.

Now let $x=0$ and $v\geq0$. Then $\lambda v\geq0$ and $\inf_{y< 0}g(y)=g(0)=0$, $\min_{y\geq0}g(y)=g(\lambda v)=1-(\lambda/2)v^2\geq0\Leftrightarrow |v|\leq\sqrt{2/\lambda}$. Hence $\levp f(0)=[0,\sqrt{2/\lambda}]$.

(ii) Resolvents are immediate due to~\ref{subdifferentials3}, so we only justifies formula for $\prox_\lambda f$. Note that $\prox_\lambda f (x)=\argmin_{y\in\R}g(y)$, where $g(y)=f(y)+(1/2\lambda)(y-x)^2$ satisfies
\begin{align*}
g(y)=
\begin{cases}
\frac{1}{2\lambda}(y-x)^2,&\text{if }y\leq0;\\
\frac{1}{2\lambda}(y-x)^2+1,&\text{if }y>0.
\end{cases}
\end{align*}
When $x\leq0$, it is easy to see that $\argmin_{y\in\R}g(y)=\{x\}$. Now consider $x>0$, in which case $g(0)=\min_{y\leq0}g(y)=(1/2\lambda)x^2$ and $g(x)=\min_{y>0}g(y)=1$. Then $g(0)<g(x)\Leftrightarrow x<\sqrt{2\lambda}$ and $g(0)=g(x)\Leftrightarrow x=\sqrt{2\lambda}$, from which the desired formula readily follow. \endproof

\begin{example}\label{-l2 example} Let $\lambda>0$ and $f(x)=-\norm{x}$ for $x\in\Rn$. Then the following hold:
	\begin{enumerate}
		\item\label{negative l2 subdiff}
		$
		\partial f(x)=
		\begin{cases}
	\sph, &\text{if }x=0;\\
-x/\norm{x}, &\text{if }x\neq0,
		\end{cases}
		$,
		$
		\partial_pf(x)=
		\begin{cases}
		\emptyset,&\text{if }x=0;\\
		-x/\norm{x}, &\text{if }x\neq0,
		\end{cases}
		$, and
		$$
\levp f(x)=
		\begin{cases}
			\emptyset,&\text{if }\norm{x}<\lambda;\\
		-x/\norm{x}, &\text{if }\norm{x}\geq\lambda.
		\end{cases}
		$$
		\item \label{negative l2 resolvent}
		$
\prox_\lambda f(x)
		=\left(\Id+\lambda\levp f\right)^{-1}(x)
		\subset\left(\Id+\lambda\partial_p f\right)^{-1}(x)
		\subset\left(\Id+\lambda\partial f\right)^{-1}(x)
		$,
		where
		$$
	\left(\Id+\lambda\partial_p f\right)^{-1}(x)=
	\begin{cases}
\{x\pm\lambda\proj_\sph(x)\},&\text{if }\norm{x}<\lambda,\\
	\{x+\lambda\proj_\sph(x)\},&\text{if }\norm{x}\geq\lambda,
	\end{cases}
		$$

		$$
		\left(\Id+\lambda\partial f\right)^{-1}(x)=
		\begin{cases}
\{x\pm\lambda\proj_\sph(x)\},&\text{if }\norm{x}\leq\lambda,\\
			\{x+\lambda\proj_\sph(x)\},&\text{if }\norm{x}>\lambda.
		\end{cases}
		$$
		and
		$$(\forall x\in\Rn)~
		\prox_\lambda f(x)=x+\lambda\proj_\sph(x).
		$$
\item\label{i:negative-norm} $f$ is $1/\lambda$-hypoconvex on convex subsets of $\{x\in\Rn:\ \|x\|\geq \lambda\}$.
	\end{enumerate}
\end{example}
\proof \ref{negative l2 subdiff} Formulas for $\partial f$ and $\partial_p f$ follow easily from definition, hence we only justify the one for $\levp f$.
Due to Proposition~\ref{prop: basic propeties}\ref{relation2}, it suffices to check whether $(\forall x\neq0)$ $v=-x/\norm{x}\in\levp f(x)$ holds.
In turn,
$v\in\levp f(x)\Leftrightarrow
(\forall y\in\Rn)~-\norm{y}\geq-\norm{x}=\ip{v}{y-x}-(1/2\lambda)\norm{y-x}^2\Leftrightarrow
\min_{y\in\Rn}g(y)\geq0,
$
where
$
g(y)=(1/2\lambda)\norm{y-x}^2+\ip{v}{x-y}-\norm{y}+\norm{x}.
$

When $0<\norm{x}<\lambda$,
$g(x+2\lambda v)=
\norm{x}-\norm{x+2\lambda v}
=\norm{x}-|\norm{x}-2\lambda|
=2\norm{x}-2\lambda<0$,
implying that $\levp f(x)=\emptyset$.

Now let $\norm{x}\geq\lambda$.
Clearly $g$ is coercive, which means that there exists a minimizer $y\in\Rn$ of $g$.
We claim that $g(y)\geq0$, entailing the desired inclusion $v\in\levp f(x)$.
Indeed, if $y=0$, then $g(y)=(1/2\lambda)\norm{x}^2\geq0$.
If $y\neq0$, then optimality condition yields that
$
0\in\partial g(y)
=
(1/\lambda)(y-x)-v+\partial(-\norm{\cdot})(y)
$,
implying
$
y=x+\lambda(v-d),
$
where $
d
=
-y/\norm{y}
=\partial(-\norm{\cdot})(y)
$.
Hence $g(y)
=
\norm{x}-\norm{x+\lambda(v-d)}$.
Recall that $d=-y/\norm{y}$ and $v=-x/\norm{x}$.
Then
\begin{align*}
\ip{v}{d}\leq1
&\Rightarrow
\lambda-\norm{x}\leq(\lambda-\norm{x})\ip{v}{d}
\Leftrightarrow
2\lambda^2
-2\lambda^2\ip{v}{d}
+
2\lambda\ip{x}{v-d}\leq0\\
&\Leftrightarrow
\norm{x+\lambda(v-d)}^2
=
\norm{x}^2
+
\norm{\lambda(v-d)}^2
+
2\lambda\ip{x}{v-d}
\leq\norm{x}^2\\
&\Leftrightarrow
g(y)\geq0,
\end{align*}
as desired.

\ref{negative l2 resolvent} Beginning from showing $\prox_\lambda f(x)=x+\lambda\proj_\sph(x)$, we note that
\begin{equation}
(\forall y\in\Rn)~
-\norm{y}
+\frac{1}{2\lambda}\norm{y-x}^2
\geq
\frac{1}{2\lambda}\norm{y}^2
-
\left(
\frac{1}{\lambda}\norm{x}+1
\right)\norm{y}
+\frac{1}{2\lambda}\norm{x}^2\geq
-\norm{x}-\frac{\lambda}{2}.
\end{equation}
Then $(\forall y\in\Rn)~(\forall d\in\proj_\sph(x))$
\begin{equation}\label{negative l2:prox optimal value}
	-\norm{x+\lambda d}
	+
	\frac{1}{2\lambda}\norm{x+\lambda d-x}^2
	=-\norm{x}-\frac{\lambda}{2}
	\leq
	-\norm{y}
	+
\frac{1}{2\lambda}\norm{y-x}^2,
\end{equation}
implying
$x+\lambda\proj_\sph(x)
\subseteq
\prox_\lambda f(x).
$
To justify the converse, take $v\in\prox_\lambda f(x)$.
Then optimality condition implies that
\begin{equation}\label{negative l2: 1}
	\frac{1}{\lambda}(x-v)\in\partial(-\norm{\cdot})(v)\Rightarrow\norm{v-x}=\lambda.
\end{equation}
If $x=0$, then $\norm{v}=\lambda$ as desired. Now suppose $x\neq0$.
We have learned from~(\ref{negative l2:prox optimal value}) that
\begin{equation}\label{negative l2: 2}
-\norm{v}
+\frac{1}{2\lambda}\norm{v-x}^2
=
-\norm{x}
-\frac{\lambda}{2}
\Rightarrow
\norm{v}-\norm{x}
=\lambda.
\end{equation}
Combing identities~(\ref{negative l2: 1})--(\ref{negative l2: 2}) and equality characterization of the triangle inequality suggests that there exists $c\in\R$ such that $v=cx$ and
\[
|c|-1
=
|c-1|
=
\frac{\lambda}{\norm{x}}
\Rightarrow
c
=
1
+
\frac{\lambda}{\norm{x}},
\]
implying $v=x+\lambda\proj_\sph(x)$.

Next we justify $\left(\Id+\lambda\levp f\right)^{-1}=\Id+\lambda\proj_\sph$. Invoking the formula for $\levp f$ yields
$
v\in\left(\Id+\lambda\levp f\right)^{-1}(x)
\Leftrightarrow
(1/\lambda)(x-v)\in\levp f(v)
\Leftrightarrow
\norm{v}\geq\lambda$
 and $
x=(1-\lambda/\norm{v}) v.
$
We claim that
\begin{equation}\label{negative l2: equivalency}
	\norm{v}\geq\lambda
	\text{ and }
	x=\left(1-\frac{\lambda}{\norm{v}}\right)v
	\Leftrightarrow
	v\in x+\lambda\proj_\sph(x),
\end{equation}
which entails the desired formula. To establish~(\ref{negative l2: equivalency}), we argue by cases.
Suppose first the left hand side holds.
If $\norm{v}=\lambda$, then $x=0$ and clearly $v\in\lambda\proj_\sph(x)$.
If $\norm{v}>\lambda$, then $1-\lambda/\norm{v}>0$, which implies that $x/\norm{x}=v/\norm{v}$ and
$x=v-\lambda v/\norm{v}=v-\lambda x/\norm{x}$, justifying the right hand side of~(\ref{negative l2: equivalency}).
Now let $v\in x+\lambda\proj_\sph(x)$.  If $x=0$ then clearly $\norm{v}=\lambda$ and $x=(1-\lambda/\norm{v})v=0$. If $x\neq0$, then $v=x+\lambda x/\norm{x}$ and consequently $\norm{v}=\norm{x}+\lambda>\lambda$, thus completes the proof of~(\ref{negative l2: equivalency}).  Formulas for $(\Id+\lambda\partial f)^{-1}$ and $(\Id+\lambda\partial_p f)^{-1}$ can be proved similarly by repeating the above argument with extra care for the case $\norm{x}=\lambda$, thus omitted for simplicity.

\ref{i:negative-norm} Apply Proposition~\ref{p:hypoconvex}.
\endproof

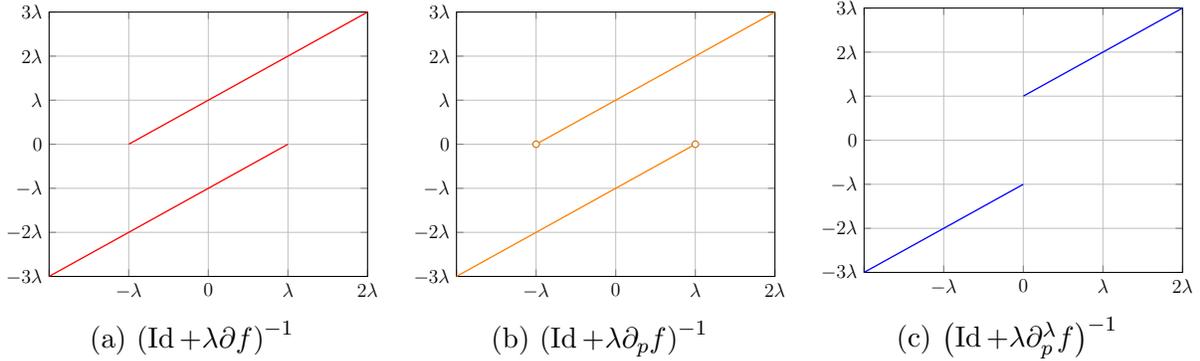
\begin{figure}[h!]
	\begin{subfigure}[b]{0.32\textwidth}
		\centering
		\resizebox{\linewidth}{!}{
			\begin{tikzpicture}
				\begin{axis}[xmin = -2, xmax = 2,ymin = -3, ymax = 3,grid = both,
					legend style={
						legend columns=2,
						anchor=south,
						at={(0.5,1)}
					},
					xtick={-2,-1,0,1,2},
					ytick={-3,-2,-1,0,1,2,3},
					xticklabels={, $-\lambda$,0,$\lambda$, $2\lambda$},
					yticklabels={$-3\lambda$,$-2\lambda$, $-\lambda$,0,$\lambda$, $2\lambda$,$3\lambda$}
					]
					\addplot[domain =-1:2,samples = 10,smooth,thick,red] {x+1};
					\addplot[domain = -2:1,samples = 10,smooth,thick,red] {x-1};
				\end{axis}
			\end{tikzpicture}
		}
		\caption{$\left(\Id+\lambda\partial f\right)^{-1}$}
		\label{fig:subfig8}
	\end{subfigure}
	\begin{subfigure}[b]{0.32\textwidth}
		\centering
		\resizebox{\linewidth}{!}{
			\begin{tikzpicture}
				\begin{axis}[xmin = -2, xmax = 2,ymin = -3, ymax = 3,grid = both,
					legend style={
						legend columns=2,
						anchor=south,
						at={(0.5,1)}
					},
					xtick={-2,-1,0,1,2},
					ytick={-3,-2,-1,0,1,2,3},
					xticklabels={, $-\lambda$,0,$\lambda$, $2\lambda$},
					yticklabels={$-3\lambda$,$-2\lambda$, $-\lambda$,0,$\lambda$, $2\lambda$,$3\lambda$}
					]
					
					\addplot[domain =-1:2,samples = 10,smooth,thick,orange] {x+1};
					\addplot[domain = -2:1,samples = 10,smooth,thick,orange] {x-1};
					\addplot[color=orange!80!black, only marks, style={mark=*, fill=white}] coordinates {(1,0)};
					\addplot[color=orange!80!black, only marks, style={mark=*, fill=white}] coordinates {(-1,0)};
				\end{axis}
			\end{tikzpicture}
		}
		\caption{$\left(\Id+\lambda\partial_p f\right)^{-1}$}
		\label{fig:subfig9}
	\end{subfigure}
	\begin{subfigure}[b]{0.32\textwidth}
		\centering
		\resizebox{\linewidth}{!}{
			\begin{tikzpicture}
				\begin{axis}[xmin = -2, xmax = 2,ymin = -3, ymax = 3,grid = both,
					legend style={
						legend columns=2,
						anchor=south,
						at={(0.5,1)}
					},
					xtick={-2,-1,0,1,2},
					ytick={-3,-2,-1,0,1,2,3},
					xticklabels={, $-\lambda$,0,$\lambda$, $2\lambda$},
					yticklabels={$-3\lambda$,$-2\lambda$, $-\lambda$,0,$\lambda$, $2\lambda$,$3\lambda$}]
					
					\addplot[domain =0:2,samples = 10,smooth,thick,blue] {x+1};
					\addplot[domain = -2:0,samples = 10,smooth,thick,blue] {x-1};
				\end{axis}
			\end{tikzpicture}
		}
		\caption{$\left(\Id+\lambda\partial_p^\lambda f\right)^{-1}$}
		\label{fig:subfig10}
	\end{subfigure}
	\caption{Resolvents in Example~\ref{-l2 example} with $n=1$.}
\end{figure}

\begin{remark} For a general Lipschitz function $f:\R\rightarrow \R$, note that the set $\dom\partial_{p}f$ can be a countable set \cite{bgw98}, and that $\partial f(x)=[-1,1]$ for every $x\in\R$ \cite{wang04}. In such a pathological case,
$(\Id+\partial f)^{-1}(x)=x+[-1,1]$ for every $x\in \R$, which is certainly not the proximal mapping $P_{1}f$.
\end{remark}


\section*{Acknowledgments}
Xianfu Wang and Ziyuan Wang were supported by NSERC Discovery grants.

\section*{Data availability statement}
All data generated or analysed during this study are included in this article.

\bibliographystyle{siam}
\bibliography{prox}

\begin{thebibliography}{10}

\bibitem{BC}
{\sc H.~H. Bauschke and P.~L. Combettes}, {\em Convex Analysis and Monotone
  Operator Theory in {H}ilbert Spaces}, Springer, Cham, 2017.

\bibitem{bmw21}
{\sc H.~H. Bauschke, W.~Moursi, and X.~Wang}, {\em Generalized monotone
  operators and their averaged resolvents}, Mathematical Programming, 189
  (2021), pp.~55--74.

\bibitem{beck2017first}
{\sc A.~Beck}, {\em {First-order Methods in Optimization}}, SIAM, 2017.

\bibitem{bh96}
{\sc J.~Benoist and J.-B. Hiriart-Urruty}, {\em What is the subdifferential of
  the closed convex hull of a function?}, SIAM Journal on Mathematical
  Analysis, 27 (1996), pp.~1661--1679.

\bibitem{bt05}
{\sc F.~Bernard and L.~Thibault}, {\em Prox-regular functions in hilbert
  spaces}, Journal of Mathematical Analysis and Applications, 303 (2005),
  pp.~1--14.

\bibitem{bgw98}
{\sc J.~M. Borwein, R.~Girgensohn, and X.~Wang}, {\em On the construction of
  {H\"older} and proximal subderivatives}, Canadian Mathematical Bulletin, 41
  (1998), pp.~497--507.

\bibitem{cwp20}
{\sc J.~Chen, X.~Wang, and C.~Planiden}, {\em A proximal average for
  prox-bounded functions}, SIAM Journal on Optimization, 30 (2020),
  pp.~1366--1390.

\bibitem{clarke98}
{\sc F.~H. Clarke, Y.~S. Ledyaev, R.~J. Stern, and P.~R. Wolenski}, {\em
  Nonsmooth Analysis and Control Theory}, Springer-Verlag, New York, 1998.

\bibitem{clarke_proximal_1995}
{\sc F.~H. Clarke, R.~J. Stern, and P.~R. Wolenski}, {\em Proximal smoothness
  and the lower-{$C^2$} property}, Journal of Convex Analysis, 2 (1995),
  pp.~117--144.

\bibitem{mor2018variational}
{\sc B.~S. Mordukhovich}, {\em Variational Analysis and Applications},
  Springer, Cham, 2018.

\bibitem{rockafellar_variational_1998}
{\sc R.~T. Rockafellar and R.~J.~B. Wets}, {\em Variational {Analysis}},
  vol.~317, Springer, Berlin.

\bibitem{wang04}
{\sc X.~Wang}, {\em Subdifferentiability of real functions}, Real Analysis
  Exchange, 30 (2004/05), pp.~137--171.

\bibitem{wang2022mirror}
{\sc Z.~Wang, A.~Themelis, H.~Ou, and X.~Wang}, {\em A mirror inertial
  forward-reflected-backward splitting: {global convergence and linesearch
  extension beyond convexity and Lipschitz smoothness}}, arXiv:2212.01504,
  (2022).

\end{thebibliography}

\end{document}